\documentclass[a4paper,11pt,makeidx]{amsart}
\usepackage{amscd}
\usepackage{xypic}  
\usepackage{amssymb}
\usepackage{amsthm}
\usepackage{epsf}
\makeindex
\newtheorem{thm}{Theorem}[section]
\newtheorem{cor}[thm]{Corollary}
\newtheorem{example}[thm]{Example}
\newtheorem{lemma}[thm]{Lemma}
\newtheorem{prop}[thm]{Proposition}
\newtheorem{defn}[thm]{Definition} 

\newtheorem{rem}[thm]{Remark}

\newcommand{\eprf}{\hfill $\square$ \bigskip\par}

\def\Sing{\mbox{Sing }}

\def\rk{\mbox{rk }}

\def\c1{\mbox{c_1 }}
\def\c2{\mbox{c_2 }}

\def\GG{\mathbb G}

\def\PP{{\mathbb P}}
\def\ZZ{{\mathbb Z}}
\def\QQ{{\mathbb Q}}
\def\CC{{\mathbb C}}
\def\O{{\mathcal O}}

\def\E{{\mathcal E}}
\def\H{{\mathcal H}}
\def\G{{\mathcal G}}

\def\cH{{\mathcal H}}
\def\cO{{\mathcal O}}
\def\L{{\mathcal L}}
\def\V{{\mathcal V}}

\def\N{{\mathcal N}}
\def\T{{\mathcal T}}
\def\FF{{\mathbb F}}
\def\E{{\mathcal E}}
\def\G{{\mathcal G}}
\def\x{\times}                   
                 
\def\iso{\simeq}                 
\def\eqv{\equiv}
\def\sub{\subseteq}

\def\+{\oplus}                  
\def\*{\otimes}                  
       
\def\law{\longrightarrow}

\def\coker{\mbox{coker }}

\def\tC0{\tilde{C_0}}
\def\tA1{\tilde{A_1}}
\def\tA2{\tilde{A_2}}

\def\Pic{\mbox{Pic }}
\def\Num{\mbox{Num }}
\def\Div{\mbox{Div }}

\def\deg{\mbox{deg }}

\begin{document}

\title{On varieties which are uniruled by lines}

\author{A.~L.~Knutsen, C.~Novelli, A.~Sarti}

\address{\hskip -.43cm Andreas Leopold Knutsen, Department of Mathematics, University of Oslo, PO
Box 1053, Blindern NO-0316 Oslo, Norway. e-mail {\tt andreas.knutsen@math.uio.no}}

\address{\hskip -.43cm Carla Novelli, Dipartimento di Matematica, Universit{\`a} degli Studi di Trento,
via Sommarive, 14, I-38050 Povo (TN), Italy. e-mail {\tt novelli@science.unitn.it}}

\address{\hskip -.43cm Alessandra Sarti, FB 17 Mathematik und Informatik, Universit\"at Mainz,
Staudingerweg 9, 55099 Mainz, Germany. e-mail {\tt sarti@mathematik.uni-mainz.de}}

\maketitle

\begin{abstract}
 Using the $\sharp$-minimal model program of uniruled varieties we show that
 for any pair $(X, \H)$ consisting of a reduced and irreducible variety $X$ of dimension $k \geq 3$ and a globally generated big line bundle $\H$ on $X$ with $d:= \H^k$ and $n:= h^0(X, \H)-1$ such that $d<2(n-k)-4$, then $X$ is uniruled of $\H$-degree one, except if $(k,d,n)=(3,27,19)$ and a ${\sharp}$-minimal model of $(X, \H)$ is $(\PP^3,\O_{\PP^3}(3))$. We also show that the bound is optimal for threefolds.
\end{abstract}

\renewcommand{\thefootnote}{\alph{footnote}}

\footnote[0]{ 14H45.}

\section{Introduction}

It is well-known that an irreducible nondegenerate variety $X \sub \PP^n$ of degree $d$ satisfies $d \geq n - \dim X+1$. Varieties for which equality is obtained are the well-known {\it varieties of minimal degree}, which are completely classified.

Varieties for which $d$ is ``small'' compared to $n$ have been the objects of intensive study throughout the years, see e.g. \cite{ha,ba,F1,F2,F3, is, io, ho,reid, mella}. One of the common features is that such varieties are covered by rational curves.

More generally one can study pairs $(X,\H)$ where $X$ is an irreducible variety (possibly with some additional assumptions on its singularities) and $\H$  a line bundle on $X$ which is sufficiently ``positive'' (e.g. ample or (birationally) very ample or big and nef). Naturally we set $d:= \H^{\dim X}$ and $n:= \dim |\H|$. The difference between $d$ and $n$ is measured by the $\Delta$-genus: $\Delta(X,\H):=d+\dim X-n-1$, introduced by Fujita (cf. \cite{F1} and 
\cite{F2}) and Fujita in facts shows that $\Delta(X,\H) \geq 0$ for $X$ smooth and $\H$ ample.
The cases with $\Delta(X,\H) =0$ are the varieties of minimal degree and the cases with 
$\Delta(X,\H) =1$ have been classified by Fujita \cite{F1,F2,F3} and Iskovskih \cite{is}.

The notion of being covered by rational curves is incorporated in the concept of a variety being {\it uniruled}. Roughly speaking a variety is uniruled if through any point there passes a rational curve. In many ways these are the natural generalizations to higher dimensions of ruled surfaces.
In the Mori program they play an important role, because - like in the case of ruled surfaces - these are the varieties for which the program does not yield a minimal model, but a Mori fiber space. Uniruled varieties can also be considered to be the natural generalizations to higher dimensions of surfaces of negative Kodaira dimension: in fact it is conjectured that a (smooth) variety is uniruled if and only if its Kodaira dimension is negative. The conjecture has been established for threefolds by Miyaoka \cite{mi}.

With the evolution of a structure theory for higher dimensional varieties in the past decades, namely the Mori program, the geometry of rational curves on varieties has gained new importance. The main idea is to obtain information about varieties by studying the rational curves on them (cf. e.g. \cite{kol}). 

To measure the ``degree'' of the rational curves which cover $X$ we say in addition that $X$ is uniruled of $\H$-degree at most $m$ if the covering curves all satisfy $\Gamma \cdot \H \leq m$. 

For surfaces Reid \cite{reid} found bounds on the uniruledness degree of $(X, \H)$ depending on $d$ and $n$. For instance he shows that $X$ is uniruled of $\H$-degree one if $d < \frac{4}{3}n-3$. The same result was obtained by Horowitz \cite{ho} using a different approach. In particular, it immediately follows (by taking surface sections) that if $\dim X=3$ and $|\H|$ contains a smooth surface then $X$ is uniruled of $\H$-degree one if 
$d < \frac{4}{3}(n-1)-3=\frac{4}{3}n-\frac{13}{3}$, and similarly for higher 
dimensions. (Note that there is even a better bound one can obtain from the results on threefolds in \cite{ho}, see Proposition \ref{horobound} below). However it is to be expected that this ``naive'' inductive procedure does not yield an optimal bound. 

The purpose of this article is to obtain a bound for uniruledness degree one that is optimal for threefolds. In fact we show:

\begin{thm} \label{main}
Let $(X, \H)$ be a pair consisting of a reduced and irreducible $3$-dimensional variety $X$ and a globally generated big 
line bundle $\H$ on $X$. Set $d:= \H^3$ and $n:= h^0(X, \H)-1$.

If $d<2n-10$ then $X$ is uniruled of $\H$-degree one, except when $(d,n)=(27,19)$ and a ${\sharp}$-minimal model of $(X, \H)$ is $(\PP^3,\O_{\PP^3}(3))$.
\end{thm}

(For the definition of a ${\sharp}$-minimal model we refer to Definition \ref{def:MM} below.)

The bound in Theorem \ref{main} is sharp since there are pairs satisfying $d=2n-10$ for infinitely
many $d$ and $n$, namely $(\PP^2 \times \PP^1, \O_{\PP^2}(2) \boxtimes
\O_{\PP^1}(a))$ (cf. Example \ref{sharp} below), which are not uniruled of $\H$-degree one. 

Observe that in Corollary
\ref{cor:easybound} below we obtain a better bound than in Theorem \ref{main}
for $n \leq 12$. 

As  a consequence of Theorem \ref{main} we get the following result for higher dimensional varieties, which is probably far from being sharp:

\begin{cor} \label{maingen}
Let $(X, \H)$ be a pair consisting of a reduced and irreducible $k$-dimensional variety $X$, $k \geq 4$,  and a globally generated big line bundle $\H$ on $X$. Set $d:= \H^k$ and $n:= h^0(X, \H)-1$. If  $d<2(n-k)-4$ then $X$ is uniruled of $\H$-degree one.
\end{cor}

For those preferring the notion of $\Delta$-genus, the condition $d  <2(n-k)-4$ is equivalent to $\Delta(X, \H) < n-k-5= h^0(X, \H)- \dim X - 6$. 

The above results have the following corollary for embedded varieties:

\begin{cor} \label{mainproj}
Let $X \subset {\mathbb P}^n$ be a nondegenerate irreducible
variety of dimension $k \geq 3$ and degree $d$. If $d<2(n-k)-4$ then $X$ is uniruled by lines, except when $(k,d,n)=(3,27,19)$ and a ${\sharp}$-minimal model of $(X, \O_X(1))$ is 
$(\PP^3,\O_{\PP^3}(3))$.
\end{cor}

Note that the condition $d<2(n-k)-4$ implicitly requires $n \geq k+6$ in Corollary \ref{mainproj}.

To prove these results we use the $\sharp$-minimal model program of uniruled
varieties introduced for surfaces by Reid in \cite{reid} and developed for threefolds by Mella in \cite{mella}. The main advantage of the $\sharp$-minimal model program is that one does not only work with birational modifications along the minimal model program but also 
uses a polarizing divisor. Under certain assumptions one manages to follow every step of the program on an effective divisor, i.e. a (smooth) surface in the case of threefolds.

Our method of proof uses the classification results in \cite{mella} and borrows ideas from \cite{reid}. The crucial point is a careful investigation of pairs $(X, \H)$
such that the output of the $\sharp$-minimal model program is a particular type of Mori fiber space which we call a {\it terminal Veronese fibration} (see Definition
\ref{VFSdefi} below): this  is roughly speaking a terminal threefold marked by a
line bundle with at most base points fibered
over a smooth curve with general fibers being smooth Veronese surfaces (with respect to the marking line bundle) and
having at most finitely many fibers being cones over a smooth quartic
curve. We find a lower bound on the degree of such a
threefold (in fact on every marked terminal threefold having a terminal Veronese fibration as a
$\sharp$-minimal model) and on the number of degenerate fibers of the members of 
the marking linear system. 

The precise statement, which we hope might be of independent interest, is the following:

\begin{prop} \label{VFS}
Let $(X, \H)$ be a three-dimensional terminal Veronese fibration 
(see Definition \ref{VFSdefi}) over a smooth curve $B$ and set 
$n:=h^0(\H)-1$ and $d:= \H^3$.
Then $d \geq 2n-10$ and the general member of $|\H|$
is a smooth surface fibered over $B$ with $ \geq  \frac{n-5}{2}$ fibers which are unions of two 
conics (with respect to $\H$) intersecting in one point (the other
fibers are smooth quartics).
\end{prop}

Observe that both equalities are obtained by 
$(\PP^2 \times \PP^1, \O_{\PP^2}(2) \boxtimes\O_{\PP^1}(a))$, cf. Examples \ref{sharp} and \ref{comesharp}. 

In Section \ref{basic} we set notation and give all central definitions. In particular we introduce, after \cite{mella}, the ${\sharp}$-minimal models of uniruled threefolds in Subsection \ref{subbasic} and give a result that will be central for our purposes in Lemma \ref{lemma:eximini}.

In Section \ref{biguni} we study threefolds with big uniruled systems, as defined in \cite{mella}, and after obtaining some basic results (such as Lemmas \ref{cliffcon} and \ref{corner2}) we first obtain an ``easy bound'' on $d$ such that a threefold is uniruled 
in degree one (Proposition \ref{easybound}) and then we show how to reduce the proofs of our main results
Theorem \ref{main} and its two corollaries to a result about uniruled threefolds having a terminal Veronese fibration as a
$\sharp$-minimal model, namely Proposition \ref{boundX}.

The proofs of Proposition \ref{boundX} and of Proposition \ref{VFS} are then settled in Section \ref{hyperplane}. 

Finally, in Section \ref{final} we mention a couple of slight improvements of some results in \cite{mella} that follow from our treatment (Propositions 
\ref{impro1} and \ref{impro2}). \\

\noindent{\bf Acknowledgments} We are indebted to M.~Mella for suggesting the problem and for many helpful discussions. We were introduced to the topic during the wonderful
school {\it Pragmatic 2002} in Catania, and it is a great pleasure to thank all the participants, as well as the organizer A.~Ragusa. 

Moreover we thank M.~Andreatta, R.~Pignatelli, W.~Barth and K.~Ranestad 
for useful comments.

\section{Basic definitions and notation} \label{basic}

We work over the field of complex numbers. 

\subsection{Uniruled varieties} An $n$-dimensional projective variety $X$ is called {\it uniruled} if there is a variety $Y$ of dimension $n-1$ and a generically finite dominant rational map $p: Y \x \PP^1 -- \rightarrow X$. In particular, such a variety is covered by rational curves (cf. \cite[IV 1.4.4]{kol}). 

Note that since our field is algebraically closed the condition is equivalent to the map $p$ being a birational morphism (cf. \cite[IV (1.1.1)]{kol}).

If $H$ is a nef $\QQ$-divisor on $X$ and $d \in \QQ$ we say that $X$ is {\it uniruled of $H$-degree at most $d$} if $\deg (p^*H)_{| \PP^1 \x \{y\} } \leq d$ for every $y \in Y$, or equivalently if there is a dense open subset $U \sub X$ such that every point in $U$ is contained in a rational curve $C$ with $C.H \leq d$ (cf. \cite[IV 1.4.]{kol}). A consequence is that in fact {\it every}
point in $X$ is contained in a rational curve $C$ with $C.H \leq d$ (cf. \cite[IV (1.4.4)]{kol}).
In particular, if $X \sub \PP^n$ we say that $X$ is uniruled by lines if $d=1$ with respect to $H:=\O_X(1)$

We say that a linear system $\mathfrak{D}$ of Weil divisors on $X$ is  {\it movable} if for the general member $H \in \mathfrak{D}$ we have $h^0(X, nH) >1$ for some $n >0$. We say that it is {\it nef} if the general member $H \in \mathfrak{D}$ is nef.
Assume furthermore that $X$ is terminal $\QQ$-factorial. Then the 
{\it threshold} of the pair $(T,\mathfrak{D})$  is defined as
\[
\rho_{\mathfrak{D}}=\rho(T,\mathfrak{D}):=\mbox{sup}\left\{m \in \mathbb{Q} ~ : ~H+mK_T \mbox{~is~an~effective}~\mathbb{Q}\mbox{-divisor}\right\}\geq 0
\]

For a pair $(X, \L)$ of a variety $X$ and a line bundle $\L$ on $X$ we similarly define
\[ \rho_{\L}=\rho(T,\L):= \rho(T,|\L|). \]
Moreover we set
\begin{equation}\label{eq:ned}
 d(X,\L):=\L^{\dim X} \;  \mbox{and} \; n(X,\L):=h^0(X,\L)-1= \dim |\L|. 
\end{equation}
In these terms the $\Delta$-genus, introduced by Fujita (cf. \cite{F1} and 
\cite{F2}), is 
\begin{equation}\label{eq:deltagenus}
 \Delta(X, \L):= d(X,\L)+ \dim X -n(X,\L)-1 
\end{equation}
and all the results in the paper can be equivalently formulated with the $\Delta$-genus.

\subsection{${\sharp}$-Minimal models of uniruled threefolds} \label{subbasic}

Recall that a surjective morphism $f: X \law Y$ with connected fibers between normal varieties is called a {\it Mori fiber space} if $-K_X$ is $f$-ample, $\rk \Pic (X/Y)=1$ and $\dim X > \dim Y$. Given $A,B\in \Div(X)\otimes \QQ$ then A is {\it f-numerically equivalent} to B ($A \eqv_f B$) if $A\cdot
C=B\cdot C$ for every curve contracted by f.
The following result of Mella will be central to us:

\begin{thm} \label{eximini}{\rm \cite[Thm. 3.2]{mella}} 
Let $T$ be a terminal $\QQ$-factorial uniruled threefold and $\mathfrak{D}$ a movable and nef linear system. 
Then there exists a pair $(T^{{\sharp}}, \mathfrak{D} ^{{\sharp}})$, called a ${\sharp}$-minimal model of the pair $(T, \mathfrak{D})$, with general element $S ^{{\sharp}} \in \mathfrak{D} ^{{\sharp}}$ such that:
  \begin{itemize}
  \item[(i)]   There is a birational map $\phi: T --\rightarrow T^{{\sharp}}$ with $\phi_* \mathfrak{D} = \mathfrak{D} ^{{\sharp}}$.
  \item[(ii)] $\pi: T^{{\sharp}} \law W$ is a Mori fiber space.
  \item[(iii)] $\rho(T,\mathfrak{D})K_{T^{{\sharp}}}+ S ^{{\sharp}} \eqv _{\pi} \O_{T^{{\sharp}}}$.
  \end{itemize}
\end{thm}

\begin{rem} \label{rem:eximini}
 {\rm  Note that it is implicitly shown in the proof of \cite[Thm. 3.2]{mella} that
$\rho(T,\mathfrak{D})= \rho(T^{{\sharp}},\mathfrak{D}^{{\sharp}})$.}
\end{rem}

The following explains the terminology used in Theorem \ref{main} and Corollaries 
\ref{maingen} and \ref{mainproj}:

\begin{defn} \label{def:MM}
  {\rm For any pair $(X, \L)$ consisting of a reduced and irreducible uniruled variety $X$ and a big and nef line bundle $\L$ on $X$ we will by a {\it ${\sharp}$-minimal model of $(X, \L)$} mean a  
${\sharp}$-minimal model of $(\tilde{X}, |f^*\L|)$, where $f: \tilde{X} \law X$ is a minimal resolution of singularities.}
\end{defn}

As observed in \cite[Prop. 3.6]{mella}, if $\rho(T,\mathfrak{D}) <1$ and there is a smooth surface $S \in \mathfrak{D}$, one can describe each step in the ${\sharp}$-minimal model program in a neighborhood of $S$. More precisely, 
set $T_0:=T$, $\mathfrak{D}_0:=\mathfrak{D}$, $S_0:=S$ and $T_m:=T^{\sharp}$, $\mathfrak{D}_m:=\mathfrak{D}^{{\sharp}}$ and $S_m=S^{{\sharp}}:= \phi_*S$. Denote by $\phi_i: T_{i-1} --\rightarrow T_{i}$ for $i=1, \ldots, {m}$ each birational modification in the ${\sharp}$-minimal model program relative to $(T,\mathfrak{D})$ and define inductively $S_i:= \phi_*S_{i-1}$ and $\mathfrak{D}_i:= \phi_*\mathfrak{D}_{i-1}$. Then each $S_i$ is smooth and the map 
\[ {\phi_i}_{|S_{i-1}}: S_{i-1} \law S_i \]
is a well-defined morphism which is either an isomorphism or a contraction of a disjoint union of 
$(-1)$-curves. 
We will need the following precise result, especially the inequality in (\ref{eq:nu1}), which in principle is the same as \cite[Claim 2.6]{reid} and which might have some interest of its own:

\begin{lemma} \label{lemma:eximini}
With the same notation and assumptions as in Theorem \ref{eximini}, assume furthermore that
$\rho:=\rho(T,\mathfrak{D}) <1$ and there is a smooth surface $S \in \mathfrak{D}$. Set $S^{{\sharp}}:= \phi_*S$. 

Then both $S$ and $S^{\sharp}$ are smooth and Cartier and
the map 
\[ f:=\phi_{|S}: S \law S^{\sharp} \]
is a well-defined morphism which can be factorized in a finite sequence of monoidal transformations. 

Let $l$ be the total number of irreducible curves contracted by $f$.

If $\rho \geq \frac{1}{3}$ then, setting $D:= \O_S(S)$, we have
\begin{equation}
\label{eq:nu1}
(D+\frac{\rho}{1-\rho}K_S)^2 \geq -l(\frac{\rho}{1-\rho})^2. 
\end{equation}
\end{lemma}

\begin{proof}
By what we discussed above it is clear that $f$ is a well-defined morphism which can be factorized into a  succession of $l$ monoidal transformations and that $S^{\sharp}$ is smooth. Moreover, since $S^{{\sharp}}$ is smooth and $T^{\sharp}$ is terminal, we have
$S^{{\sharp}} \cap \Sing T^{\sharp}= \emptyset$ \cite[2.3]{mella}, so that
$S^{{\sharp}}$ is Cartier.
Setting $f_i:={\phi_i}_{|S_i}$ we can factorize $f$ as
\[
\xymatrix{
 S \ar[r]^{f_1} & S_1 \ar[r]^{f_2} &  \cdots \ar[r]^{f_{m-1}} & S_{m-1} \ar[r]^{f_m} & S_m=  S^{{\sharp}},
}
\] 
where each $f_i$ contracts $l_i$ disjoint $(-1)$-curves $E^i_1, \ldots, E^i_{l_i}$ with $l_i \geq 0$. 
Clearly $l = \sum_{i=1}^m l_i$.

We set $D_i:= \O_{S_i}(S_i)$ and $D^{{\sharp}}:= \O_{S^{\sharp}}(S^{\sharp})$. 

If $\phi_i$ is a flip then $S_i$ is disjoint from the flipping curves by \cite[Claim 3.7]{mella}, so
that $f_i$ is an isomorphism.

If $\phi_i$ contracts a divisor onto a curve then it is shown in \cite[Case 3.8]{mella} that 
the fiber $F_i$ of $\phi_i$ satisfies $S_i \cdot F_i=0$, whence $D_i \cdot F_i=0$, which means that
all $E^i_j$ satisfy $E^i_j \cdot D_i=0$. 

If $\phi_i$ is a divisorial contraction onto a curve then it is shown in \cite[Case 3.9]{mella} that 
$f_i$ is a contraction of a single $(-1)$-curve $E_i=E^i_1$ which satisfies $E_i \cdot D_i=1$. 

In other words, for every $j$ we have three possibilities:
\begin{itemize}
\item[(a)] $l_i =0$.
\item[(b)] $l_i >0$ and $E^i_j \cdot D_i=0$ for all $j$.
\item[(c)] $l_i =1$ and $E^i_1 \cdot D_i=1$.
\end{itemize}

Now denote by $L^i_j$ the total transform of $E^i_j$ on $S$. Then $(L^i_j)^2=-1$ and 
$L^i_j \cdot L^{i'}_{j'}=0$ for $(i,j) \neq (i',j')$. We have
\begin{equation}
  \label{eq:nu11} 
K_S = f^*K_{S^{\sharp}} + \sum L^i_j 
\end{equation}
and by (a)-(c)   
\begin{equation}
  \label{eq:nu10} 
D = f^*D^{\sharp} - \sum \mu^i_j L^i_j \; \mbox{with} \; \mu^i_j \in \{0,1 \}.\end{equation}
Hence
\[
D+\frac{\rho}{1-\rho}K_S = f^*(D^{\sharp}+ \frac{\rho}{1-\rho}K_{S^{\sharp}}) +
    \sum (\frac{\rho}{1-\rho}-\mu^i_j) L^i_j, 
\]
and since there are $l$ terms in the sum we get
\begin{equation}
  \label{eq:nu2}
 (D+\frac{\rho}{1-\rho}K_S)^2 = (D^{\sharp}+ \frac{\rho}{1-\rho}K_{S^{\sharp}})^2
                           -l(\frac{\rho}{1-\rho})^2 + \sum \mu^i_j(\frac{2\rho}{1-\rho} - \mu^i_j).
\end{equation}

By definition and invariance of $\rho$ (cf. Remark \ref{rem:eximini}) we have that 
$\rho K_{T^{{\sharp}}}+ S ^{{\sharp}}$ is $\QQ$-effective. From Theorem \ref{eximini}(iii) we have that it is also $\QQ$-nef, whence its restriction to $S^{\sharp}$ is also $\QQ$-effective and $\QQ$-nef.
Since $S^{{\sharp}}$ is Cartier  we get by
adjunction that $(\rho K_{T^{{\sharp}}}+ S ^{{\sharp}})_{|S^{\sharp}} \iso 
(1-\rho)(\frac{\rho}{1-\rho}K_{S^{{\sharp}}}+D^{\sharp})$, whence by $\QQ$-nefness
\begin{equation}
  \label{eq:nu3} 
 (D^{\sharp}+ \frac{\rho}{1-\rho}K_{S^{\sharp}})^2 \geq 0. 
\end{equation}
Moreover the assumption $\rho \geq \frac{1}{3}$ is equivalent to
$\frac{\rho}{1-\rho} \geq \frac{1}{2}$, whence
\begin{equation}
  \label{eq:nu4} 
 \sum \mu^i_j(\frac{2\rho}{1-\rho} - \mu^i_j) \geq \mu^i_j(1-\mu^i_j) \geq 0.
\end{equation}
Now (\ref{eq:nu1}) follows combining (\ref{eq:nu2})-(\ref{eq:nu4}) 
\end{proof}

A direct consequence is the following

\begin{cor}{\rm \cite[Cor. 3.10]{mella}} \label{cor:eximini}
With the same notation and assumptions as in Theorem \ref{eximini}, assume furthermore that
$\rho:=\rho(T,\mathfrak{D}) <1$ and that $\mathfrak{D}$ is base point free. 
Then $\mathfrak{D}^{\sharp}$ is a linear system of Cartier divisors with at most base points and with general element a smooth surface.
\end{cor}

The following result will be useful to us:

\begin{lemma}\label{laccadegree}
With the same notation and assumptions as in Thm 1.1, assume furthermore
that $\rho(T, \mathfrak{D}) < 1$ and there is a smooth surface $S \in \mathfrak{D}$.
Let $S^{{\sharp}}:= \phi_*S$.

If $T^{\sharp}$ is uniruled of $S^{\sharp}$-degree at most $m$
then also  $T$ is uniruled of $S$-degree at most $m$.
\end{lemma}

\begin{proof}
By assumption $T^{\sharp}$ is covered by a family of rational curves
$\Gamma$ such that $\Gamma \cdot S^{\sharp} \leq m$. The strict transform
$\tilde{\Gamma}$ of each such $\Gamma$ on $T$ then satisfies
$\tilde{\Gamma} \cdot S \leq m$ by \cite[Lemma 3.15]{mella}.
Moreover, as we saw in the proof of Lemma \ref{lemma:eximini}, the
exceptional locus of $\phi$ is covered by rational curves
$\Delta$ satisfying $\Delta \cdot S = 0$ or $1$ and the result follows.
\end{proof}

\section{Threefolds with big uniruled systems} \label{biguni}

Let $T$ be a terminal $\QQ$-factorial threefold and $\mathfrak{D}$ a movable linear system on $T$. Recall from
\cite[Def. 5.1]{mella} that $(T, \mathfrak{D})$ is said to be a pair with a big uniruled system if there is an $S \in \mathfrak{D}$ which is big and nef and a smooth surface of negative Kodaira dimension.
Such a pair is uniruled with $\rho(T,\mathfrak{D}) <1$ by \cite[Lemma 5.2]{mella}.

Note that by Lemma \ref{lemma:eximini} $S$ is Cartier. Therefore $\mathfrak{D}$ is a linear system of Cartier divisors and $\H:=\O_T(S)$ is a line bundle. Moreover $\mathfrak{D}$ is automatically movable 
since $S$ is big and nef. 
We will say that 
$(T, \H)$ is a pair with a big uniruled system if $(T, |\H|)$ is. Moreover, since the general 
$S^{\sharp} \in \mathfrak{D}^{\sharp}$ is Cartier by Lemma \ref{lemma:eximini} 
$\H^{\sharp}:=\O_{T^{\sharp}}(S^{\sharp})$ is a line bundle and we will say that the pair
$(T^\sharp,{\H}^{\sharp})$ is a ${\sharp}$-minimal model of the pair $(T, \H)$.

\begin{thm} \label{me1} {\rm \cite[Thm. 5.3]{mella}}
Let $T$ be a terminal $\QQ$-factorial threefold and $\H$ a line bundle on $T$ such that
$(T,\H)$ is a pair with a big uniruled system. Set $\rho:= \rho(T,\H)$.

Then $(T^\sharp, \H^{\sharp})$ is one of the following:

\begin{itemize}
\item[{\rm (I)}] a $\mathbb{Q}$-Fano threefold of index $1/\rho>1$, with 
$K_{T^\sharp} \sim -1/\rho \H^{\sharp}$, belonging to Table \ref{tabella}. 
\item[{\rm (II)}] a bundle over a smooth curve with generic fiber $(F,\H^{\sharp}_{|F})\cong(\mathbb{P}^2,\cO(2))$ and with at most finitely many fibers $(G,\H^{\sharp}_{|G})\cong (\bf{S}_4,\cO(1))$, where $\bf{S}_4$ is the cone over the normal quartic curve and the vertex sits over a hyperquotient singularity of type $1/2(1,-1,1)$ with $f=xy-z^2+t^k$, for $k\geq 1$. ($\rho=2/3$)
\item[{\rm (III)}] a quadric bundle with at most $c A_1$ singularities of type $f=x^2+y^2+z^2+t^k$, for $k\geq 2$, and $\H_{|F}^{\sharp}\sim\cO (1)$. ($\rho=1/2$)
\item[{\rm (IV)}] $(\mathbb{P}(\E),\cO(1))$ where $\E$ is a rank $3$ vector bundle over a smooth curve. ($\rho=1/3$)
\item[{\rm (V)}] $(\mathbb{P}(\E), \cO(1))$ where $\E$ is a rank $2$ vector bundle over a surface of negative Kodaira dimension. ($\rho=1/2$)
\end{itemize}
\end{thm}

In the cases (I) the general $S ^ {{\sharp}} \in |\H ^ {{\sharp}}|$ is a smooth del Pezzo surface and $\O_{S^ {\sharp}}(\H ^ {{\sharp}}) \iso \frac{\rho}{\rho-1}K_{S^ {\sharp}}$. A list of such threefolds (with corresponding values for $\rho$) is given in \cite{cafl}. Moreover one can easily calculate $d(T^{{\sharp}}, \H ^ {\sharp})$  and $n(T^{{\sharp}}, \H ^ {\sharp})$. Indeed 
\[ d(T^{{\sharp}}, \H ^ {\sharp}) := (\H ^ {\sharp})^3= (\O_{S ^ {\sharp}}(\H ^ {{\sharp}}))^2 = \frac{\rho^2}{(\rho-1)^2} K_{S^ {\sharp}}^2,\]
and by Riemann-Roch
\[ n(T^{{\sharp}}, \H ^ {\sharp}):= h^0(\H ^ {\sharp})-1= h^0(\O_{S^ {\sharp}}(\H ^ {{\sharp}}))= 
\frac{\rho}{2(\rho-1)^2} K_{S^ {\sharp}}^2+1.\] 

In Table \ref{tabella} we list all the cases (see \cite[p.~81]{cafl}). In the table  ${\mathbb P}(w_1,\ldots, w_n)$ denotes the weighted projective space with weight $w_i$ at the coordinate $x_i$. The hyperplane given by $x_i$ is denoted 
$\{x_i=0\}$. Moreover $X_a$ (resp. $X_{a,b}$) denotes a hypersurface of degree $a$
(resp. a complete intersection of two hypersurfaces of degrees $a$ and $b$) and similarly for $H_a$ and $H_{a,b}$. 
The variety $\GG(1,4)$ is the Grassmannian parameterizing lines in $\PP^4$, embedded in 
${\mathbb P}^9$ by the Pl\"ucker embedding.

\begin{footnotesize}
\begin{table}\caption{$\QQ$-Fano threefolds}\label{tabella}
\begin{tabular}{|c|c|c|c|c|c|c|c|c|} 
\hline 
Type &  $T^ {\sharp}$ & general $S^{{\sharp}} \in |\H^ {{\sharp}}|$ & $\rho$   & $-\frac{\rho}{\rho-1}$ & 
$K_{H^ {\sharp}}^2$ & $d(T^{{\sharp}}, \H ^ {\sharp})$  & $n(T^{{\sharp}}, \H ^ {\sharp})$  
\\ \hline  \hline
(a) &  $X={\mathbb P}(1,1,2,3)$  & $H_6 \subset X$ & $6/7$  & $6$ & $1$ & 36 & 22 \\ \hline 
(b) & $X_6 \subset {\mathbb P}(1,1,2,3,3)$ & $X_6 \cap \{x_4=0\})$ & $3/4$ & $3$ & $1$ & $9$ & $7$ \\
 \hline
(c) & $X_6 \subset {\mathbb P}(1,1,2,3,4)$ & $X_6 \cap \{x_4=0\}$ & $4/5$ & $4$ & $1$ & $16$ & $11$ \\
 \hline
(d) & $X_6 \subset {\mathbb P}(1,1,2,3,5)$ & $X_6 \cap \{x_4=0\}$ & $5/6$ & $5$ & $1$ & $25$ & $16$ \\
 \hline
(e) & $X_6\subset {\mathbb P}(1,1,2,2,3)$ & $X_6 \cap \{x_3=0\}$ & $2/3$ & $2$ & $1$ & $4$ & $4$ \\
 \hline
(f) & $X_6\subset {\mathbb P}(1,1,1,2,3)$ & $X_6 \cap \{x_0=0\}$ & $1/2$ & $1$ & $1$ & $1$ & $2$ \\
 \hline
(g) & $X={\mathbb P}(1,1,1,2)$ & $H_4 \subset X$ & $4/5$ & $4$ & $2$ & $32$ & $21$ \\
\hline
(h) & $X_4 \subset {\mathbb P}(1,1,1,2,2)$ & $X_4\cap \{x_4=0\}$ & $2/3$ & $2$ & $2$ & $8$ & $7$ \\
 \hline
(i) & $X_4 \subset {\mathbb P}(1,1,1,2,3)$ & $X_4\cap \{x_4=0\}$ & $3/4$ & $3$ & $2$ & $18$ & $13$ \\
 \hline
(j) & $X_4 \subset {\mathbb P}(1,1,1,1,2)$ & $X_4 \cap \{x_0=0\}$ & $1/2$ & $1$ & $2$ & $2$ & $3$ \\
\hline
(k) & $X={\mathbb P}^3$ & $H_3 \subset X$ & $3/4$ & $3$ & $3$ & $27$ & $19$ \\
\hline
(l) & $X_3 \subset {\mathbb P}(1,1,1,1,2)$ & $X_3 \cap \{x_4=0\}$ & $2/3$ & $2$ & $3$ & $12$ & $10$ \\
 \hline
(m) & $X_3 \subset{\mathbb P}^4$ & $X_3 \cap \{x_0=0\}$ & $1/2$ & $1$ &  $3$ & $3$ & $4$ \\
 \hline
(n) & $X_2 \subset{\mathbb P}^4$ & $H_{2,2} \subset X_2$ & $2/3$ & $2$ & $4$ & $16$ & $13$ \\
 \hline
(o) & $X_{2,2} \subset{\mathbb P}^5$ &  $X_{2,2} \cap \{x_0=0\}$ & $1/2$ & $1$ &  $4$ & $4$ & $5$ \\   \hline
(p) & $X={\mathbb P}^6 \cap \GG (1,4) $ & $X \cap \{x_0=0\}$ & $1/2$ & $1$ &  $5$ & $5$ & $6$ \\ 
\hline
(q) & $X_2 \subset{\mathbb P}^4$ & $X_2 \cap \{x_0=0\} \iso {\mathbb P}^1 \x {\mathbb P}^1$ & $1/3$ & $1/2$ & $8$ & $2$ & $4$ \\ \hline
(r) & $X={\mathbb P}^3$ & ${\mathbb P}^1 \x {\mathbb P}^1 \iso H_2 \subset X$ & $1/2$ & $1$ & $8$ & $8$ & $9$ \\ \hline
(s) & $X={\mathbb P}^3$ & $ \{x_0=0\} \iso {\mathbb P}^2 \subset X$ & $1/4$ & $1/3$ & $9$ & $1$ & $3$ 
\\ \hline
(t) & $X={\mathbb P}(1,1,1,2)$ & $\{x_3=0\} \iso {\mathbb P}^2 \subset X$ & $2/5$ & $2/3$ & $9$ & $4$ & $6$ \\ \hline
\end{tabular}
\end{table}
\end{footnotesize}

The following easy consequence of Clifford's theorem will be useful for our purposes:

\begin{lemma}\label{cliffcon}
Let $(X, \H)$ be a pair consisting of an irreducible terminal $\QQ$-factorial threefold $X$ 
and a globally generated and big line bundle $\H$ on $X$. Set $d:= d(X,\H)$ and $n:=n(X,\H)$.

If $d < 2n-4$ then $(X, \H)$ is a pair with a big uniruled system and for any smooth irreducible surface $S \in |\H|$ and any irreducible curve $D \in |\H_{|S}|$ we have
\begin{equation} \label{eq:p5} 
D \cdot K_S \leq d-2n+2.
\end{equation}
\end{lemma}

\begin{proof}
The general element $S \in |\H|$ is a smooth irreducible surface.

Pick any irreducible curve $D \in |\O_S(\H)|$. 
Then $\deg \O_D(\H)= \H^3=d$ and
from
\begin{equation}
  \label{eq:p0}
0 \law \O_X  \law \H \law  \O_S(\H) \law 0  
\end{equation}
and
\begin{equation}
  \label{eq:p1}
0 \law \O_S  \law \O_S(\H) \law  \O_D(\H) \law 0  
\end{equation}
we get 
\begin{equation}
  \label{eq:p2}
 h^0(\O_D(\H)) \geq h^0(\O_S(\H))-1 \geq h^0(\H)-2=n-1. 
\end{equation}
Hence
\begin{equation}
  \label{eq:p3}
 \deg \O_D(\H) -2(h^0(\O_D(\H))-1) \leq  d-2(n-2) < 0,
\end{equation}
whence by Clifford's theorem on irreducible singular curves (see the appendix of \cite{eks}) we must have 
$h^1(\O_D(\H))=0$, so that $\chi(\O_D(\H))=h^0(\O_D(\H)) \geq  n-1$. 
We see from (\ref{eq:p1}) that $h^1(\O_S(\H)) \leq h^1(\O_S)$, so 
that we have, from (\ref{eq:p2}),
\[
  \chi (\O_S(D))- \chi (\O_S) =   \chi (\O_{D}(\H)) =h^0(\O_{D}(\H)) \geq   n-1.
\]
Combining with Riemann-Roch we get
\begin{eqnarray*}
 D \cdot K_S & = & D^2-2(\chi (\O_S(D))- \chi (\O_S)) \\
     &  \leq  & d-2n+2 < -2,
\end{eqnarray*}
proving (\ref{eq:p5}) and showing that $\kappa(S) <0$, so that $(X, \H)$ is a pair with a big uniruled system.
\end{proof}

\begin{lemma}\label{corner2}
With the same assumptions and notation as in Lemma \ref{cliffcon}, assume furthermore that 
$(X^{{\sharp}}, \H ^ {\sharp})$ is of type (I) in Theorem \ref{me1}. Then, setting 
$\rho:=\rho(X, \H)$, we have, for general $S^{\sharp} \in |\H ^ {\sharp}|$ (which is a smooth del Pezzo surface),
\begin{equation}
  \label{eq:p52}
  d-2n+2= \frac{\rho}{\rho-1}K_{S^ {\sharp}}^2 + \sum \mu^i_j,
\end{equation}
where the $\mu^i_j$ are defined as in (\ref{eq:nu10}).
\end{lemma}

\begin{proof}
  Pick smooth surfaces $S \in |\H|$ and $S^{\sharp} \in |\H^{\sharp}|$ such that we have a 
commutative diagram
\[
\xymatrix{
 X \ar@{-->}[r]^{\phi} &  X^{\sharp} \\ 
 S \ar[r]^f \ar@{^{(}->}[u]  & S^{\sharp} \ar@{^{(}->}[u]
}
\]
where $f:= \phi_{|S}$ is well-defined and is a successive contraction of 
disjoint unions of finitely many $(-1)$-curves (in each step) by Lemma \ref{lemma:eximini}. Since $S^{\sharp}$ is a smooth del Pezzo surface, we have $h^1(\O_S)=h^1(\O_{S^{\sharp}})=0$. It is then easily seen by the proof of Lemma \ref{cliffcon} that equality holds in 
(\ref{eq:p5}) (note that we have $h^1(\O_X) \leq h^1(\O_X(-\H))+h^1(\O_S)=0$ by Kawamata-Viehweg vanishing). Using (\ref{eq:nu11}) and (\ref{eq:nu10}) we therefore get, for $D \in |\O_S(\H)|$:
\begin{eqnarray*}
  d-2n+2 & = & D.K_S = (f^*D^{\sharp} - \sum \mu^i_j L^i_j).(f^*K_{S^{\sharp}} + \sum L^i_j) 
   \\  
         & = & D^{\sharp}.K_{S^{\sharp}}+\sum \mu^i_j = \frac{\rho}{\rho-1}K_{S^ {\sharp}}^2 + \sum \mu^i_j,
\end{eqnarray*}
as stated. 
\end{proof}

As a ``warming up'' before proceeding with the proofs of the main results  we give the proof of the following bound. Note that it improves Theorem \ref{main} for $n \leq 12$.

\begin{prop} \label{easybound}
Let $(X, \H)$ be a pair consisting of an irreducible terminal $\QQ$-factorial threefold $X$ 
and a globally generated and big 
line bundle $\H$ on $X$. Set $d:= d(X,\H)$ and $n:=n(X,\H)$.

If $n \geq 4$, $d<\frac{4}{3}n-\frac{4}{3}$, $d \neq n-1$ for $n \leq 9$ and $d \neq n-2$ for $n \leq 6$, then $X$ is uniruled of $\H$-degree one.
\end{prop}

\begin{proof}
Since $n \geq 4$ we have $2n-4 \geq \frac{4}{3}n-\frac{4}{3}$, whence $d <
2n-4$ so by Lemma \ref{cliffcon} $(X, \H)$ is a pair with a big uniruled system. 

Moreover for general $S \in |\H|$ and $D \in |\H_{|S}|$ we have, by Lemma \ref{cliffcon},
\begin{eqnarray*}
  D \cdot \Big(\frac{3}{2} \H+K_X\Big) & = & D \cdot (\H+K_X) + \frac{1}{2}D \cdot \H = 
D \cdot K_S + \frac{1}{2}d \\
 & \leq & d-2n+2 +  \frac{1}{2}d = \frac{3}{2}d -2n +2 \\
   & < & \frac{3}{2} \Big( \frac{4}{3}n-\frac{4}{3} \Big)-2n +2 = 0, 
\end{eqnarray*}
whence $\rho(X,\H) < 2/3$.

It follows that the ${\sharp}$-minimal model $(X^\sharp,{\cH}^\sharp)$ is in the list of 
Theorem \ref{me1} and moreover it cannot be as in (II) since $\rho(X,\cH) = 2/3$ in this case.
In the cases (III)-(V) one immediately sees that $(X^\sharp,{\cH}^\sharp)$ is uniruled of $\H ^\sharp$-degree one, whence $(X,{\cH})$ is also uniruled of $\H$-degree one by Lemma \ref{laccadegree}.

We have $n(T^{{\sharp}}, \H ^ {\sharp}) \geq n \geq 4$, and by using Table \ref{tabella}
we see that the cases in (I) where $n(X^\sharp,{\cH}^\sharp) \geq 4$ and $\rho < 2/3$ are the cases (m), (o), (p), (q), (r) and (t). Among these (m), (o), (p) and (q) are clearly uniruled of $\H ^\sharp$-degree one. This leaves us with the cases (r) and (t).

By (\ref{eq:p52}) and Table \ref{tabella} we have 
$d-2n= \sum \mu^i_j -10$ in case (r) and $d-2n= \sum \mu^i_j -8$ in case (t).
By using (\ref{eq:nu10}) we obtain $d=d(T^{{\sharp}}, \H ^ {\sharp})- \sum (\mu^i_j)^2=d(T^{{\sharp}}, \H ^ {\sharp})- \sum \mu^i_j$ since $\mu^i_j\in\{0,1\}$ and we find
that $n= 9- \sum \mu^i_j$ in case (r) and 
$n=6-\sum \mu^i_j$ in case (t). It follows that $d=n-1$ for $n \leq 9$ in case (r) and 
$d=n-2$ for $n \leq 6$ in case (t). 
\end{proof}

We have the following corollary:

\begin{cor} \label{cor:easybound}
 Let $(X, \H)$ be a pair consisting of a reduced and irreducible $3$-dimensional variety
$X$ and a globally generated and big 
line bundle $\H$ on $X$. Set $d:= d(X,\H)$ and $n:=n(X,\H)$. 

If $d<\frac{4}{3}n-\frac{4}{3}$ and $d \neq n$ when $5 \leq n \leq 8$, 
$d \neq n-1$ for $n \leq 9$ and $d \neq n-2$ for $n \leq 6$,then $X$ is uniruled by lines.
\end{cor}

\begin{proof}
Let $\pi: T \law X$ be a resolution of the singularities of $X$.
Then $\pi^*\mathcal{H}$ is globally generated and big with 
$d(T,\pi^*\mathcal{H})=\H^3=d$ and $n(T,\pi^*\mathcal{H})= \dim |\pi^*\mathcal{H}| \geq n$ and
we can apply Proposition \ref{easybound}. The additional cases $d=n$ for $5 \leq n \leq 8$
occur since
equality does not need to occur in $n(T,\pi^*\mathcal{H}) \geq n$.
\end{proof}

It is interesting to note that similar bounds as in Proposition
\ref{easybound} can be obtained by using results of Horowitz \cite{ho} and
Reid \cite{reid} respectively. As an example we briefly show the bound one
would get using Horowitz's result. Note that 
 the obtained bound is stronger than Proposition \ref{easybound} for $n\geq 29$ but one has to assume that the
threefold is smooth and the line bundle is very ample.

\begin{prop} \label{horobound}
  Let $X$ be a smooth threefold and $\H$ a very ample line bundle on $X$ with
  $d:=\H^3$ and $n:=h^0(\H)-1$. If 
$n \geq 12$ and $d<\frac{3}{2}(n-4)$, or $7 \leq n \leq 11$ and
  $d<\frac{4}{3}(n-3)$ then $X$ is uniruled of $\H$-degree one.
\end{prop}

\begin{proof}
By \cite[Corollary p.~668]{ho}, if $d<\frac{3}{2}(n-4)$ then $X$ is ruled by planes or quadrics and if $d<\frac{4}{3}(n-3)$ then $X$ is ruled by planes, in which cases it is clearly uniruled by lines.
Now we note that $\frac{4}{3}(n-3) > \frac{3}{2}(n-4)$ for $n \leq 11$ and
that the condition $d \geq n-2$ requires $n \geq 7$.  
\end{proof}

Now we give the main ideas and the strategy of the proof of Theorem \ref{main}. The main result we will need to prove is the following:

\begin{prop} \label{boundX}
Let $(X, \H)$ be a pair consisting of an irreducible terminal $\QQ$-factorial threefold $X$ 
and a globally generated, big 
line bundle $\H$ on $X$. Set $d:= d(X,\H)$ and $n:=n(X,\H)$.

If a $\sharp$-minimal model of $(X, \H)$ is of type (II) in Theorem \ref{me1}, then
$d \geq 2n-10$.
\end{prop}

The proof of this result will be given in  Section \ref{hyperplane} below, after a careful study of the threefolds of type (II) in Theorem \ref{me1}. We will now give the proofs of Theorem \ref{main} and Corollaries \ref{maingen} and \ref{mainproj} assuming 
Proposition \ref{boundX}.
\vspace{0.5cm}

\noindent \textbf{Proof of Theorem \ref{main}.}
Let $X$ be a reduced and  irreducible $3$-dimensional variety and $\H$ a globally generated big 
line bundle on $X$. Set $d:= \H^3$ and $n:= h^0(X, \H)-1$ and assume
$d<2n-10$.

Let $\pi: T \law X$ be a resolution of the singularities of $X$.
Then $\pi^*\mathcal{H}$ is globally generated and big with 
$d(T,\pi^*\mathcal{H})=\H^3=d$ and $n(T,\pi^*\mathcal{H})= \dim |\pi^*\mathcal{H}| \geq n$.
Since $(d,n)=(27,19)$ satisfies $d=2n-11$ we can reduce to the case where $X$ is smooth.
Therefore we henceforth assume $X$ is smooth.

By Lemma \ref{cliffcon} $(X, \H)$ is a pair with a big uniruled system, whence any ${\sharp}$-Minimal model 
$(X^ {{\sharp}}, \H^ {{\sharp}})$ of $(X, \H)$ is in the list of Theorem \ref{me1}. Moreover, by Proposition \ref{boundX}, it cannot be of type (II).

We easily see that the cases (III)-(V) are uniruled of $\H ^\sharp$-degree one. In the cases 
(I) we have, by (\ref{eq:p52}),
\[ \frac{\rho}{\rho-1}K_{S^ {\sharp}}^2 + \sum \mu^i_j =     d-2n+2 < -8. \]
By checking Table \ref{tabella} one finds that we can only be in case (k) with 
$\sum \mu^i_j=0$. Hence $d=d(T^{{\sharp}}, \H ^ {\sharp})=27$ by (\ref{eq:nu10}) and
$n=19$ by (\ref{eq:p52}). Now the result follows from Lemma \ref{laccadegree}.
\eprf

\noindent \textbf{Proof of Corollary \ref{maingen}.}
Let $X$ be a reduced and  irreducible variety  of dimension $k \geq 4$ and $\H$ a globally generated big  
line bundle on $X$ with $d:= \H^k$ and $n:= h^0(X, \H)-1$.

As just mentioned in the proof of Theorem \ref{main} we can assume $X$ is smooth.

Setting $X_k:=X$ and $\H_k:=\H$, we recursively choose general smooth ``hyperplane sections'' 
$X_{i-1} \in |\H_{i}|$ and define $\H_{i-1}:=\H_i \* \O_{X_{i-1}}$, for $2 \leq i \leq k$. (Note that $\dim X_i=i$ and $\H_i$ is a line bundle on $X_i$.) 

Let $n_3:=h^0(\H_3)-1$. Then from the exact sequence
\begin{equation}
  \label{eq:ma1}
  0 \law \O_{X_i} \law \H_i \law \H_{i-1} \law 0
\end{equation}
we have 
\begin{equation}
  \label{eq:ma2}
  n_3 \geq n-(k-3)= n-k+3.
\end{equation}
Together with the conditions $d < 2(n-k)-4$ this implies $d < 2n_3-10$ and it follows from
Theorem \ref{main} that either $(X_3, \H_3)$ is uniruled of  degree one or 
$(d,n_3)=(27,19)$ and $(X_3 ^{\sharp}, \H_3 ^{\sharp})$ is $({\mathbb P}^3,\cO (3))$.

In the first case, since $X_3$ is the general $3$-fold section, it follows that 
$(X, \H)$ is also uniruled of  degree one.

In the second case we have equality in (\ref{eq:ma2}), i.e.
\begin{equation}\label{eq:ma3} 
19=n_3=h^0(\H_3)-1=h^0(\H)-(k-3)-1=n-k+3.
\end{equation} 
Denote by $f:X_3 --\rightarrow X_3^{\sharp}$ the birational map of the
$\sharp$-minimal model program. By Lemma \ref{lemma:eximini} its restriction $g$ to $X_2$ is a morphism onto a smooth surface $X_2^{\sharp} \in |\H_3^{\sharp}|$ which is either an isomorphism or a finite sequence of contractions of disjoint unions of $(-1)$-curves.

Now $\H_3^{\sharp}$ is globally generated, whence the restriction 
$\H_2^{\sharp}:= \H_3^{\sharp} \* \O_{X_2^{\sharp}}$ is also globally generated. 
Moreover, we have $\dim |\H_3| = \dim |\H_3^{\sharp}|=19$, whence 
$|\H_3^{\sharp}| = f_*|\H_3|$. By \cite[Lemma 4.5]{cafl} the restriction map
$H^0(\H_3^{\sharp}) \rightarrow H^0(\H_2^{\sharp})$ is surjective and by 
(\ref{eq:ma2}) and (\ref{eq:ma3}) the restriction map
$H^0(\H_3) \rightarrow H^0(\H_2)$ is surjective as well. Hence 
$|\H_2^{\sharp}| = g_* |\H_2|$ and this can only be base point free if every curve $E$ contracted by $g$ satisfies $E.\H_2=0$. Denoting by $\varphi_{\H_2}$ and $\varphi_{\H_2^{\sharp}}$ the morphisms defined by $|\H_2|$ and $|\H_2^{\sharp}|$ respectively, this implies that $\varphi_{\H_2^{\sharp}} \circ g = \varphi_{\H_2}$, in other words the images 
$\varphi_{\H_2^{\sharp}} (X_2^{\sharp})$ and $\varphi_{\H_2} (X_2)$ are the same and
since $\H_2^{\sharp}$ is very ample and $h^0(\H_2^{\sharp}) = h^0(\H_3
^{\sharp})-1=h^0(\H_3)-1=n-k+3$, this image is a smooth nondegenerate surface
in $\PP^{n-k+3}$, call it $S$, which is birational to $X_2$. Moreover, by
(\ref{eq:ma3}), the natural map $H^0(\H) \rightarrow H^0(\H_2)$ is surjective,
so $S$  is the same as the image of the restriction of the morphism $\varphi_{\H}: X \law \PP^n$ defined by  $|\H|$ to $X_2$. Since this holds for the general surface section it follows that $\varphi_{\H}$ itself is birational onto a variety 
of dimension $k$, having $S$ as a smooth surface section.

We now apply the theorem of Zak (unpublished, cf. \cite{za}) and L'vovski (cf. \cite{lv1} and \cite{lv2}) which says (in the surface case) that a smooth nondegenerate surface $S \sub \PP^N$ satisfying 
$h^0(\N_{S/{\PP^N}}(-1)) \leq N+r$ for $r < (N-3)(N+1)$ cannot be the surface section of a nondegenerate irreducible $k$-dimensional variety in
$\PP^{N+k}$ other than a cone. Since a cone is uniruled by lines, the corollary will follow if we show that
$h^0(\N_{S/{\PP^{18}}}(-1))  \leq 20$, with $S$ being the $3$-uple embedding of a smooth 
cubic surface $S_0$ in $\PP^3$.

We argue  as in \cite[p.~160-161]{glm} to compute
$h^0(\N_{S/{\PP^{18}}}(-1))$. We give the argument for the sake of the reader.

From the Euler sequence and tangent bundle sequence
\begin{eqnarray*} 
0 \law \O_S(-1) \law \CC^{19} \* \O_S \law \T_{\PP^{18}} (-1) \* \O_S \law 0 \\
0 \law \T_S(-1) \law \T_{\PP^{18}} (-1) \* \O_S \law \N_{S/{\PP^{18}}}(-1) \law 0 
\end{eqnarray*}
we find
\begin{equation}
  \label{eq:ma5}
 h^0(\N_{S/{\PP^{18}}}(-1)) \leq h^0(\T_{\PP^{18}} (-1) \* \O_S) + h^1(\T_S(-1)) \leq 19 + 
h^1(\T_{S_0}(-3)).
\end{equation}

From the tangent bundle sequence of $S_0 \sub \PP^3$
\begin{equation}
  \label{eq:ma6}
0 \law \T_{S_0}(-3) \law \T_{\PP^3} (-3) \* \O_{S_0} \law \N_{{S_0}/{\PP^3}}(-3) \law 0  
\end{equation}
and the fact that $\N_{{S_0}/{\PP^3}}(-3) \iso \O_{S_0}$, we find
\[
 h^1(\T_{S_0}(-3)) \leq 1+h^1(\T_{\PP^3} (-3) \* \O_{S_0}).
\]
In view of (\ref{eq:ma5}) it will suffice to show that $h^1(\T_{\PP^3} (-3) \* \O_{S_0})=0$.

The Euler sequence of $S_0 \sub \PP^3$ shows that $h^1(\T_{\PP^3} (-3) \* \O_{S_0})= \coker \mu$ where 
$\mu: H^0(\O_{S_0}(1)) \* H^0(\omega_{S_0}(2)) \law H^0(\omega_{S_0}(3))$ is the multiplication map. As in \cite[p.~161]{glm} one proves the surjectivity of $\mu$ by restricting to the general smooth curve section $C \sub S_0$ and showing that the map $\nu: H^0(\O_{C}(1)) \* H^0(\omega_{S_0}(2)) \law H^0(\omega_{C}(2))$ is surjective. The latter is the composition
\begin{eqnarray*}
\nu: H^0(\O_{C}(1)) \* H^0(\omega_{S_0}(2)) & \law & H^0(\O_{C}(1)) \* H^0(\omega_{C}(1)) \\
& \law & H^0(\omega_{C}(2)).
\end{eqnarray*}
The first map is surjective since $h^1(\omega_{S_0}(1)) = h^1(\O_{S_0}(-1))=0$ and the second is surjective by Green's $H^0$-lemma \cite[(4.e.1)]{gre}.
This concludes the proof of the corollary.
\eprf

It is immediate that Corollary \ref{mainproj} follows from the two above.

As we already noted in the introduction, Theorem \ref{main} is sharp by the following example:

\begin{example}\label{sharp}
{\rm The bound of Theorem \ref{main} is sharp. In fact consider $X=\PP^2
  \times \PP^1$ with projections $p$ and $q$ respectively
  and let $\H := p^*\O_{\PP^2}(2) \otimes
q^*\O_{\PP^1}(a)$  for $a\in \mathbb{N}_{>0}$. We have
$n:=h^0(\H)-1=h^0(\mathcal{O}_{\mathbb{P}^2}(2))\cdot
h^0(\mathcal{O}_{\mathbb{P}^1}(a))-1=6(a+1)-1$ and $d:=\H^3= 
(p_2^*\mathcal{O}(2) \otimes p_1^*\mathcal{O}(a))^3 = 3(p_2^*\mathcal{O}(2))^2
  \cdot p_1^*\mathcal{O}(a)=12a$, whence $d=2n-10$.

Clearly any curve $C$ on $X$ satisfies $C \cdot \H =
C \cdot p^*\O_{\PP^2}(2) + C \cdot q^*\O_{\PP^1}(a) \geq 2$, with equality obtained 
exactly for the lines in the $\PP^2$-fibers,  so that $X$ is uniruled of $\H$-degree two and not uniruled of  $\H$-degree one.
}
\end{example}


\section{Terminal Veronese fibrations} \label{hyperplane}
  
In this section we will prove Propositions \ref{VFS} and \ref{boundX}.

Since we will have to
study the threefolds as in (II) of Theorem
\ref{me1} we find it convenient to make the following definition:

\begin{defn} \label{VFSdefi}
 {\rm Let $(X, \H)$ be a pair satisfying the following:}
\begin{itemize}
\item[(i)] {\rm $X$ is a three-dimensional terminal
 irreducible variety with a Mori fiber
 space structure $p:X \law B$, where $B$ is a smooth curve.}
\item[(ii)] {\rm $\H$ is a line bundle on $X$ such that the system $|\H|$ contains a smooth surface and has at most base points and $\H^3
 >0$.}  
\item[(iii)] {\rm The general fiber of $p$ is $(V, \H _{|V}) \iso (\mathbb{P}^2,\cO(2))$ and the rest
  are at 
  most finitely many fibers $(G, \H _{|G}) \iso (\mathbf{S}_4,\cO(1))$, where $\mathbf{S}_4$ is the cone
  over a normal quartic curve.}
\end{itemize}
{\rm Such a Mori fiber space will be called a} (three-dimensional) terminal Veronese fibration.
\end{defn}

The threefolds of Theorem
\ref{me1}~(II) are terminal Veronese fibrations. 

The easiest examples of terminal Veronese fibrations are the smooth ones in
Example \ref{sharp}. But there are also singular such varieties and these were
erroneously left out in both \cite[Prop. 3.7]{mellaad} and
\cite[Prop. 3.4]{cafl}, as remarked by Mella in \cite[Rem. 5.4]{mella}: Take
$\PP^2 \x \PP^1$ and blow up a conic $C$ in a fiber and contract the strict
transform of $C$, thus producing a Veronese cone singularity.

Although our main aim is to prove Proposition \ref{boundX} we believe that
terminal Veronese fibrations are interesting in their own rights. In order to
prove Proposition \ref{boundX} we will study ``hyperplane sections'' of $X$,
i.e. surfaces in $|\H|$, and show that the desired bound on the degree follows since
the general such surface has to have a certain number of degenerate fibers,
i.e. unions of two lines. What we first prove in this section is the following, which is part of the statement in Proposition \ref{VFS}:

\begin{prop} \label{VFSprimo}
  Let $(X, \H)$ be a three-dimensional terminal Veronese fibration and set 
$n:=h^0(\H)-1$ and $d:= \H^3$.

Then any smooth member of $|\H|$
is a surface fibered over $B$ with $k \geq \frac{n-5}{2}$ fibers which are unions of two smooth rational curves intersecting in one point (the other
fibers are smooth rational curves).
\end{prop}

\begin{proof}
  Denote by $\V$ the numerical equivalence class of the fiber. 
 Let $S \in |\H|$ be a smooth surface. Then, since $X$ is terminal, we have $S \cap \Sing X = \emptyset$ (cf. \cite[(2.2)]{mella}).

By property  (iii) any 
fiber of $S$ over $B$ is either a smooth quartic, a union of two conics intersecting in one point, or a double conic, all with respect to $\H$. Denote by $F$ the numerical equivalence class of the fiber. Then $F^2=0$.

If a fiber were a double conic, we could write $F \eqv 2F_0$ in $\Num S$. However, in this case we would get the contradiction $F_0.(F_0+K_S)=-1$, so this case does not occur.

In the case of a fiber which is a union of two conics intersecting in one point, we have $F \eqv F_1+F_2$ in $\Num S$, whence by adjunction both $F_i$ are $(-1)$-curves. Since $S$ can contain only 
finitely many such curves, the
general fiber of $S$ over $B$ is a smooth quartic and by adjunction $F.K_S = -2$.

Therefore $S$ has a finite number $k$ of degenerate fibers which are unions of two conics and since
these are all $(-1)$-curves we can blow down one of these curve in every fiber and reach a minimal model $R$ for $S$ which is a ruled surface over $B$. Let $g$ be the genus of $B$. Then
\[ k= K_R^2-K_S^2 = 8(1-g)-K_S^2 =8(1-g)-(K_X+ \H)^2 \cdot \H , \]
which only depends on the numerical equivalence class of $S$. Therefore, any smooth surface
in $|\H|$ has the same number of degenerate fibers.

Now note that the fibers of $S$ over the finitely many points of $B$ over which $X$ has singular fibers
are all smooth quartics, since $S \cap \Sing X = \emptyset$.

We now consider the birational map of $X$ to a smooth projective bundle $\tilde{X}$, as described in 
\cite[p.~ 699]{mella}.

Around one singular fiber $\bf{S}_4$ of $X$ over a point $p \in B$ this map is given by a succession of blow ups $\nu_i$ and contractions 
$\mu_i$:
\begin{equation}
  \label{eq:risol3}
\xymatrix{
           & Y_1 \ar[rd]^{\mu_1} \ar[ld]_{\nu_1} &             & \cdots \ar[ld]_{\nu_2} \ar[rd]^{\mu_{s-1}} &         & 
Y_s  \ar[ld]_{\nu_{s}} \ar[rd]^{\mu_{s}} & & \\
X=X_0 &                     & X_1 &                                    & X_{s-1}  &
                       & X_{s} = \tilde{X}       
} 
\end{equation}
where the procedure ends as soon as some $X_{s}$ has a fiber over $p\in B$ which is a smooth Veronese surface.

For every $\nu_i$ the exceptional divisor $E_i$ is either a smooth Veronese surface or a cone over the quartic curve, and the strict transform of the singular fiber $\bf{S}_4$ of $X_{i-1}$ over $p$ is 
$G_i \iso \FF_4$, the desingularization of $\bf{S}_4$. These two intersect along a smooth quartic $C_i$. Then $\mu_i$ contracts $G_i$ and $X_i$ is smooth along the exceptional locus of the contraction and contains the image $C_i'$ of $C_i$ which is a smooth quartic  with respect to $\H$.

Following $S$ throughout the procedure, we see that $S$ stays out of the exceptional locus of every $\nu_i$
and in the contraction it is mapped to a surface having $C_i'$ as fiber over $p$. 

In other words the procedure of desingularizing one singular fiber of $X$ maps every smooth surface in  
$|\H|$ to a smooth surface passing trough a unique smooth quartic over $p$. 

Doing the same procedure for all the other singular fibers of $X$ we therefore reach a smooth projective bundle $\PP(\E)$ over $B$ and under this process 
$|\H|$ is ``mapped'' to a (not necessarily complete) linear system on $\PP(\E)$ having smooth quartics over the corresponding points of $B$ as base curves. Denote the corresponding line bundle on 
 $\PP(\E)$ by $\H'$. Since we have not changed the number of degenerate fibers of any smooth surface in $|\H|$ over $B$, we see that every smooth surface in $|\H'|$ still has $k$ degenerate fibers over $B$. 
Since clearly $\dim |\H'| \geq \dim |\H|$ it is now sufficient to show that any smooth surface in $|\H'|$
has $k \geq  \frac{h^0(\H')-6}{2}$ fibers which are unions of two conics intersecting in one point. This is the contents of the following proposition.
\end{proof}

\begin{prop} \label{liscio}
  Let $f:X \iso \PP(\E) \law B$ be a three-dimensional projective bundle over a smooth curve of genus $g$. Denote by $\V$ the numerical equivalence class of the fiber and assume $\H$ is a line bundle on $X$ satisfying:
  \begin{itemize}
  \item[(i)]  $\H_{|V} \iso \O_{\PP^2}(2)$ for every fiber $V \iso \PP^2$,
  \item[(ii)] $|\H|$ is nonempty with general element a smooth surface, 
  \item[(iii)]  the only curves in the base locus of $|\H|$, if any, are smooth conics in the fibers.
  \end{itemize}
Then any smooth surface in $|\H|$
is fibered over $B$ with $k$ fibers which are unions of two lines intersecting in one point, where
\begin{equation}
  \label{eq:a1}
  k= \frac{1}{4}\H^3 \geq \frac{h^0(\H)-6}{2}.
\end{equation}
\end{prop}

\begin{proof}
  We only have to prove (\ref{eq:a1}).

Since every fiber of $X$ over $B$ is a $\PP^2$ we have $(K_X)_{|\PP^2} \iso K_{\PP^2} \iso \O_{\PP^2}(-3)$
  so we can choose a very ample line bundle $\G \in \Pic X$ such that
  \begin{equation} \label{eq:a2} 
\Num X \iso \ZZ \G \+ \ZZ \V, \; \G^2 \cdot \V=1, \; \G \cdot \V^2 = \V^3=0, 
    \end{equation}
  \begin{equation} \label{eq:a4} 
K_X \eqv b\V - 3\G, \; b \in \ZZ
  \end{equation}
and
  \begin{equation} \label{eq:a5} 
\H \eqv a\V +2\G, \; a \in \ZZ.
  \end{equation}

The general element $G \in |\G|$ is a smooth ruled surface over $B$; in particular
\[
 8(1-g)  =  K_G^2 = (K_X+ \G)^2 \cdot \G = (b\V-2\G)^2 \cdot \G =  4\G^3-4b,
\]
that is
\begin{equation}
  \label{eq:a6}
  \G^3=2(1-g)+b.
\end{equation}
Let now $S \in |\H|$ be any smooth surface. Clearly (as discussed in the proof of the previous proposition)
$K_S^2=8(1-g)-k$. We compute, using (\ref{eq:a2}) and (\ref{eq:a6}), 
\begin{eqnarray*}
  K_S^2 & = & (K_X+\H)^2 \cdot \H = ((a+b)\V-\G)^2 \cdot (a\V+2\G) \\
        & = & 2\G^3 -3a-4b = 2(2(1-g)+b)-3a-4b \\
        & = & 4(1-g)-3a-2b,
\end{eqnarray*}
so that 
\begin{equation}
  \label{eq:a6'}
  k= 4(1-g)+3a+2b.
\end{equation}
At the same time we have
\begin{eqnarray}
  \label{eq:a7'} \H^3 & = & (a\V+2\G)^3 = 12a +8 \G^3 \\
\nonumber            & = &  12a+8b+16(1-g) = 4k,
\end{eqnarray}
proving the equality in (\ref{eq:a1}). 

The inequality in (\ref{eq:a1}) we have left to prove is $\H^3 \geq 2h^0(\H)-12$. We therefore assume, to get a contradiction, that
\begin{equation}
  \label{eq:a7}
  \H^3 \leq 2h^0(\H)-13.
\end{equation}

Note that $\H$ is nef. Indeed, if $\gamma \cdot \H < 0$ for some irreducible curve $\gamma$ on $X$, then clearly $\gamma$ is contained in the base locus of $|\H|$, so that by assumption (iii) $\gamma$ has to be a smooth conic in some fiber $V \iso \PP^2$ of $f$. Hence
\[ \gamma \cdot \H = \gamma \cdot \H_{|V} = \O_{\PP^2}(2) \cdot \O_{\PP^2}(2)=4, \]
a contradiction. This proves that $\H$ is nef.

Let now $C \in |\H_{|S}|$ be general.

Since the $1$-dimensional part of the base locus of $|\H|$ can only consist of smooth conics in the fibers of $f$, we can write, on $S$, 
\[ C \sim C_0 + (f_{|S})^* \mathfrak{v} \eqv C_0 + cF, \]
for some nonnegative integer $c$, where $\mathfrak{v}$ is an effective divisor of degree
$c$ on $B$; $F$ denotes the numerical equivalence class of the fiber of $f_{|S}:S \law B$;
and $C_0$ is a reduced and irreducible curve (possibly singular). Consider
\begin{equation}
  \label{eq:a9}
0 \law \O_S(f^* \mathfrak{v}) \law \O_S(\H) \law \O_{C_0}(\H) \law 0.  
\end{equation}
From this sequence and 
\begin{equation}
  \label{eq:a8'}
  0 \law \O_X \law \O_X(\H) \law \O_S(\H) \law 0
\end{equation}
we get, using (\ref{eq:a7}): 
\begin{eqnarray} 
  \nonumber h^0(\O_{C_0}(\H)) & \geq & h^0(\O_{S}(\H)) - h^0(\O_S(f^* \mathfrak{v})) \\
 \label{eq:a12'}                   & \geq & h^0(\O_{X}(\H)) - h^0(B,\mathfrak{v})-1 \\
 \nonumber          & \geq &  \frac{1}{2}(\H^3+13)-c-2 = \frac{1}{2}\H^3-c+\frac{9}{2}.
\end{eqnarray}
Moreover $\deg \O_{C_0}(\H)= \H^2 \cdot (\H-c\V) = \H^3 -c\H^2 \cdot \V= \H^3-4c$, so that
\[ \deg \O_{C_0}(\H)-2(h^0(\O_{C_0}(\H))-1) \leq \H^3-4c-(\H^3-c+7) = -2c-7 < 0. \]
By Clifford's theorem on singular curves (see the appendix of \cite{eks}) we must therefore
have 
\begin{equation}
  \label{eq:a13}
  h^1(\O_{C_0}(\H))=0.
\end{equation}

Also note that since $X$ is a projective bundle over a smooth curve of genus $g$, we have
\begin{equation}
  \label{eq:a8}
  h^0(\O_X)=1, \; h^1(\O_X)=g, \; h^2(\O_X)=h^3(\O_X)=0.
\end{equation}

From Riemann-Roch on $S$ and the fact that 
$h^2(\O_S(f^* \mathfrak{v}))=h^0(K_S-f_{|S}^* \mathfrak{v})=0$ (since $F$ is nef with $F.(K_S-f^* \mathfrak{v})=F.(K_S-cF)=-2$) we find
\begin{eqnarray}
  \label{eq:a14}
  h^1(\O_S(f^* \mathfrak{v})) & = & 
  -\chi(\O_S(f^* \mathfrak{v}))+h^0(\O_S(f^* \mathfrak{v}))+h^2(\O_S(f^* \mathfrak{v})) \\
\nonumber   & = & -\frac{1}{2}cF \cdot (cF-K_S) +g-1 + h^0(\O_S(f^* \mathfrak{v})) \\
\nonumber   & \leq & -c +g-1 +c+1 = g.
\end{eqnarray}
Combining all (\ref{eq:a9})-(\ref{eq:a14}) we find
\begin{eqnarray}
  \label{eq:a15}
  h^1(\H) & \leq & h^1(\O_X)+h^1(\O_{S}(\H)) \\ 
\nonumber   & \leq & h^1(\O_X)+h^1(\O_S(f^* \mathfrak{v}))+h^1(\O_{C_0}(\H)) \\
\nonumber   & \leq &             g+g+0 =2g,
\end{eqnarray}
together with
\begin{equation}
  \label{eq:a16}
  h^2(\H)=h^3(\H)=0.
\end{equation}

We now compute $\chi(\H)$.

From
\[ 0 \law \T_G \law {\T_X}_{|G} \law \N_{G/X} \law 0 \]
we find, using (\ref{eq:a2}) and (\ref{eq:a6}), 
\begin{eqnarray*}
  c_2(X) \cdot \G & = & c_2(G)-K_G \cdot \N_{G/X} = 4(1-g) -(K_X+\G) \cdot \G^2 \\
         & = & 4(1-g) -(b\V-3\G) \cdot \G^2 - \G^3 = 4(1-g) -b +2 \G^3 \\
        & = & 4(1-g) -b +2 (2(1-g)+b) = 8(1-g)+b.
\end{eqnarray*}
Similarly we find, for any fiber $\PP^2 \iso V \sub X$,
\[
  c_2(X) \cdot \V =  c_2(V)-K_V \cdot \N_{V/X} = 3 -(K_X+\V) \cdot \V^2 =3.
\]
Hence
\begin{eqnarray}
  \label{eq:a19}
  c_2(X) \cdot \H & = & c_2(X) \cdot (2\G+a\V) = 2(8(1-g)+b)+3a \\
  \nonumber       & = & 16(1-g)+2b+3a.
\end{eqnarray}
Now Riemann-Roch on $X$ yields 
\begin{eqnarray*}
  \chi(\H) & = & \frac{1}{12}\Big \{ 2\H^3 - 3K_X \cdot \H^2 + K_X^2 \cdot \H +
         c_2(X) \cdot \H \Big \} + \chi(\O_X) \\
& = & \frac{1}{12}\Big \{ 2\Big (12a+8b+16(1-g)\Big ) - 3\Big (-12a-8b-24(1-g)\Big )  \\
& & + \; \Big (36(1-g)+6b+9a\Big )+ \Big (16(1-g)+2b+3a\Big ) \Big \} + 1-g \\
& = & \frac{1}{12}\Big \{ 72a+156(1-g) +48b \Big\} + 1-g \\
& = & 13(1-g) +6a + 4b +1-g = 14(1-g) +6a + 4b.   
\end{eqnarray*}
Comparing with (\ref{eq:a7'}) we see that
\[ \H^3 = 2\chi(\H)-12(1-g), \]
whence, using (\ref{eq:a15}) and (\ref{eq:a16}),
\begin{eqnarray*}
  \H^3 & = & 2\Big(h^0(\H)-h^1(\H)+h^2(\H)-h^3 (\H) \Big) -12(1-g) \\
       & \geq & 2\Big(h^0(\H)-2g\Big ) -12(1-g) = 2h^0(\H)-12+8g \geq 2h^0(\H)-12,
\end{eqnarray*}
contradicting (\ref{eq:a7}).

This shows that (\ref{eq:a7}) cannot hold, proving (\ref{eq:a1}).
\end{proof}

\begin{example}\label{comesharp}
{\rm As in Example \ref{sharp} take $X=\PP^2 \times \PP^1$ and $\H := p^*\O_{\PP^2}(2) \otimes q^*\O_{\PP^1}(a)$. Then  we have an embedding given by $|\H|$:
\begin{eqnarray*}
\begin{array}{ccc}
\PP^2\times\PP^1&\longrightarrow&\mathbb{P}^{6(a+1)-1}
\end{array}
\end{eqnarray*}
A hyperplane section of $X$ in $\mathbb{P}^{6(a+1)-1}$ has equation
\begin{eqnarray*}
\sum_{i,j=0,1,2, 0\leq k\leq a}l_{ijk}x_ix_jy_0^ky_1^{a-k}=0,
\end{eqnarray*}
where $(x_0:x_1:x_2)$ are the coordinates on $\PP^2$ and $(y_0:y_1)$ are the coordinates on $\PP^1$ and $l_{ijk}$ are coefficients. The section is degenerate on some Veronese surface $(\PP^2,\mathcal{O}(2))$ if the determinant of the matrix of the coefficients of the $x_ix_j$ is zero. This determinant is a polynomial of degree $3a$ in $y_0,y_1$, hence in general we find $3a$ distinct zeros. This means that  a generic hyperplane section has $3a=\frac{6(a+1)-1-5}{2}$ degenerates fibers, which is the smallest possible number of degenerate fibers for a terminal Veronese fibration as stated in Proposition \ref{VFSprimo}. }
\end{example}
\vspace{0.5cm}

\noindent \textbf{Proofs of Propositions \ref{VFS} and \ref{boundX}.} 
We note that by Proposition \ref{VFSprimo} the only statement left to prove in Proposition \ref{VFS} is a special case of Proposition \ref{boundX}.

As in Proposition \ref{boundX} let 
$(X, \H)$ be a pair consisting of an irreducible terminal $\QQ$-factorial threefold $X$ 
and a globally generated, big  
line bundle $\H$ on $X$, with $d:= d(X,\L)$ and $n:=n(X,\L)$.

Assume that a $\sharp$-minimal model $(X^{\sharp}, \H^{\sharp})$ is of type (II) in Theorem \ref{me1}, i.e. a terminal Veronese fibration over a smooth curve $B$ of genus $g$.

As we already used during the proof of Lemma \ref{corner2},
we can by Lemma \ref{lemma:eximini} pick smooth surfaces $S \in |\H|$ and $S^{\sharp} \in |\H^{\sharp}|$ such that we have a 
commutative diagram
\[
\xymatrix{
 X \ar@{-->}[r]^{\phi} &  X^{\sharp} \\ 
 S \ar[r]^f \ar@{^{(}->}[u]  & S^{\sharp} \ar@{^{(}->}[u]
}
\]
where $f:= \phi_{|S}$ is well-defined and is a succession of monoidal transformations.
Moreover $n^{\sharp} := \dim |\H^{\sharp}| \geq \dim |\H|=n$.
By Proposition \ref{VFSprimo}, $S^{\sharp}$ is fibered over $B$ with general fiber a smooth quartic and $k \geq \frac{n^{\sharp}-5}{2}$ fibers being a union of two rational curves intersecting in one point, which are both $(-1)$-curves. Therefore
\begin{equation}
  \label{eq:b1}
  K_{S^{\sharp}}^2 = 8(1-g) -k \leq 8(1-g) -\frac{n^{\sharp}-5}{2} \leq 8(1-g) -\frac{n-5}{2}.
\end{equation}

We want to show that $d \geq 2n-10$. Assume, to get a contradiction, that
\begin{equation}
  \label{eq:b2}
  d \leq 2n-11.
\end{equation}

Note that $\rho:=\rho(X, \H)= \frac{2}{3}$, so we can apply Lemma \ref{lemma:eximini}.
Let $l$ be the total number of irreducible curves contracted by $f$.
Then $K_S^2= K_{S^{\sharp}}^2-l$.
Pick any smooth irreducible curve $D \in |\O_S(\H)|$. Then by (\ref{eq:nu1}), (\ref{eq:p5}) and (\ref{eq:b1}) we have
\begin{eqnarray*} 
0 & \leq   & 4l+ (D+2K_S)^2= 4l + 4K_S^2 + 4K_S.D + D^2 \\
   &  \leq    & 4l+ 4\Big(8(1-g)-\frac{n-5}{2}-l \Big)+ 4\Big( d-2n+2 \Big) +d \\
 & =  & 32(1-g)-2(n-5) +4d -8n+8 +d \\
& =  & 32(1-g) +5d-10n+18 \leq 5d-10n+40 = 5(d-2n+10),
\end{eqnarray*}
contradicting  (\ref{eq:b2}).

We have therefore proved that $d \geq 2n-10$ and this finishes the proofs of
Propositions \ref{VFS} and \ref{boundX}.
\eprf

\section{Final Remarks} \label{final}

To conclude, we remark that a closer look at the proofs of of Propositions \ref{VFS} and \ref{boundX} shows that 
if we assume only $d <2n-4$ instead of (\ref{eq:b2}), we get $g=0$ as the only possibility.
This shows that:

{\it A three-dimensional terminal Veronese fibration over a smooth curve of genus $g>0$ must satisfy $d \geq 2n-4$.}

Consequently:

{\it If a pair $(X, \H)$ consisting of an irreducible terminal $\QQ$-factorial threefold $X$ and a globally generated, big line bundle $\H$ on $X$ has a ${\sharp}$-minimal model being of type (II) in Theorem \ref{me1} over a smooth curve of genus $g >0$, then $d \geq 2n-4$.} 

If now $(X^\sharp, \H^{\sharp})$ is a ${\sharp}$-minimal model of a pair $(X, \H)$ consisting of an irreducible terminal $\QQ$-factorial threefold $X$ and a globally generated big 
line bundle, we have seen that $\H^{\sharp}$ is still big and nef, so that 
$h^1(\O_X)=h^1(\O_S)=h^1(\O_{S^{\sharp}})=h^1(\O_{X^{\sharp}})$. We have seen that this is zero if $X^{\sharp}$ is of type (I) in Theorem \ref{boundX} and equal to $g$, the genus of $B$, if
$X^{\sharp}$ is of type (II) in Theorem \ref{boundX}. 

We have therefore obtained an improvement of \cite[Thm. 5.8]{mella}:

\begin{prop} \label{impro1}
  Let $(X, \H)$ be a pair consisting of an irreducible terminal $\QQ$-factorial threefold $X$ and a globally generated big 
line bundle $\H$ on $X$. Set $d:= \H^3$ and $n:= h^0(X, \H)-1$.

If $d<2n-10$ (resp. $d < 2n-4$ and $h^1(\O_X) >0$), then $(X^\sharp, \H^{\sharp})$ is 
of one of the types (i)-(iv) (resp. (ii)-(iv)) below: 

\begin{itemize}
\item[(i)] $(\PP^3,\O_{\PP^3}(3))$ (with $(d,n)=(27,19)$), 
\item[(ii)] a quadric bundle with at most $c A_1$ singularities of type $f=x^2+y^2+z^2+t^k$, for $k\geq 2$, and $\H_{|F}^{\sharp}\sim\cO (1)$,
\item[(iii)] $(\mathbb{P}(E),\cO(1))$ where $E$ is a rank $3$ vector bundle over a smooth curve,
\item[(iv)] $(\mathbb{P}(E), \cO(1))$ where $E$ is a rank $2$ vector bundle over a surface of negative Kodaira dimension. 
\end{itemize}
\end{prop}

Similarly, we have obtained the following improvement of 
\cite[Cor. 5.10]{mella}:

\begin{prop} \label{impro2}
  Let $(X, \H)$ be a pair consisting of an irreducible terminal $\QQ$-factorial variety $X$ of dimension $k \geq 4$ and an ample line bundle $\H$ on $X$. Set $d:= \H^k$ and $n:= h^0(X, \H)-1$.

If $d<2(n-k)-4$ then a $\sharp$-minimal model $(T^\sharp, \H^{\sharp})$ 
- in adjunction theory language the first reduction -
is either a projective bundle over a smooth curve with fibers $(\mathbb{P}^{k-1}, \O(1))$; a quadric bundle with at most $c A_1$ singularities with fibers 
$(\mathbb{Q}^{k-1}, \O(1))$; or $(\mathbb{P}(\E),\cO(1))$ where 
$\E$ is a rank $(k-1)$ ample vector bundle on $\mathbb{P}^2$ or a ruled surface. 
\end{prop}

We also observe that the two last results slightly improve a result of Ionescu, who in \cite[Theorem I]{io} describes smooth, connected, nondegenerate, complex projective varieties with degree which is small with respect to the codimension.

\end{document}